\documentclass[onefignum,onetabnum]{siamart171218}

\usepackage{amsmath}
\usepackage{amsfonts}
\usepackage{epsfig}
\usepackage{graphics} 

\usepackage{booktabs}
\usepackage{siunitx}
\usepackage{tabularx}
\usepackage{multirow}
\usepackage{subfig}
\usepackage{tikz}
\usepackage{pgfplots}

\newcommand{\Th}{\mathcal{T}_h}
\newcommand{\Eh}{\mathcal{E}_h}
\newcommand{\EK}{\mathcal{E}^K}
\newcommand{\Thk}{\mathcal{T}_h^\ell}
\newcommand{\oml}{\Omega^\ell}

\tikzset{immagine/.style={%
  above right, inner sep=0pt, outer sep=0pt}}

\ifpdf
\hypersetup{
  pdftitle={FETI-DP preconditioners for the Virtual Element Method on general 2D meshes},
  pdfauthor={D. Prada, S. Bertoluzza, M. Pennacchio, and M. Livesu}
}
\fi

\headers{FETI-DP for VEM}{D. Prada, S. Bertoluzza, M. Pennacchio, and M. Livesu}

\title{FETI-DP preconditioners for the Virtual Element Method on general 2D meshes}

\author{Daniele Prada\thanks{Istituto di Matematica Applicata e Tecnologie Informatiche del CNR, Via Ferrata 1, Pavia, Italy, (\email{daniele.prada@imati.cnr.it}, \email{silvia.bertoluzza@imati.cnr.it}, \email{micol.pennacchio@imati.cnr.it}).} \and Silvia Bertoluzza\footnotemark[1] \and Micol Pennacchio\footnotemark[1] \and Marco Livesu\thanks{Istituto di Matematica Applicata e Tecnologie Informatiche del CNR, Via dei Marini 6, Genova, Italy, (\email{marco.livesu@ge.imati.cnr.it}).}}

\newcommand{\roundPrecision}{2}

\begin{document}

\maketitle

\begin{abstract}
We analyze the performance of a state-of-the-art domain decomposition approach, the Finite Element Tearing and Interconnecting Dual Primal (FETI-DP) method~\cite{toselli2005}, for the efficient  solution of very large linear systems arising from elliptic problems discretized by the Virtual Element Method (VEM)~\cite{beirao2014}. We provide numerical experiments on a model linear elliptic problem with highly heterogeneous diffusion coefficients on arbitrary Voronoi meshes, which we modify by adding nodes and edges deriving from the intersection with an unrelated coarse decomposition. The experiments confirm also in this case that the FETI-DP method is numerically scalable with respect to both the problem size and number of subdomains, and its performance is robust with respect to jumps in the diffusion coefficients and shape of the mesh elements.
\end{abstract}

\section{Introduction}
Polytopic meshes allow the treatment of complex geometries, a crucial task for many applications in computational engineering and scientific computing. We consider here the problem of preconditioning the  Virtual Element Method (VEM), which can be viewed as an extension of the Finite Element Method to handle such a kind of meshes.
 In view of a possible parallel implementation of the method, we consider a state-of-the-art domain decomposition approach, the Finite Element Tearing and Interconnecting Dual Primal (FETI-DP) method. It has been proved that the FETI-DP method is still scalable when dealing with VEM, under the assumptions that the subdomains, obtained by agglomerating clusters of polygonal elements, are shape regular~\cite{Bertoluzza2017}. Such an assumption can be quite restrictive. In practice, it reduces to asking that the fine tessellation is built as a refinement of the previously existing coarse subdomain decomposition. This, of course, does not generally hold, so, in order to apply FETI-DP in a more general case, we propose to build the coarse decomposition independently from the tessellation, and modify the latter by inserting nodes and edges deriving from ``cutting'' it with the macro-edges of the subdomains. Of course, the resulting modified tessellation will possibly contain nasty elements with very small edges. Numerical tests do however show that FETI-DP is quite robust in this respect, providing  satisfactory results also in this framework.

This paper is organized as follows. A basic description of VEM is given in Section~\ref{sec:vem}. The FETI-DP method is introduced in Section~\ref{sec:feti-dp}, whereas the algorithm for partitioning $\Omega$ into subdomains and modifying the mesh is given in Section~\ref{sec:conformal}. Numerical experiments that validate the theory are presented in Section~\ref{sec:results}.


\section{The Virtual Element Method (VEM)}\label{sec:vem}
In this paper we focus on the numerical solution of the following model elliptic boundary value problem of second order discretized with VEM
\begin{equation}\label{eq:poisson}
-\nabla\cdot(\rho\nabla u) = f\text{ in }\Omega,\qquad u = 0\text{ on }\partial\Omega,
\end{equation}
with $f\in L^2(\Omega)$, where $\Omega\subset\mathbb{R}^2$ is a polygonal domain. We assume that the coefficient $\rho$ is a scalar such that for almost all $x \in \Omega$,  $\alpha \leq \rho(x) \leq M$ for two constants $M \geq \alpha >0$.
The variational formulation of equation~\eqref{eq:poisson} reads as follows: find $u \in V:=H^1_0(\Omega)$ such that
\begin{equation}\label{eq:variational}
		a(u,v) = ( f , v ) \ \forall v \in V, 
\end{equation}
with
\[
a(u,v) = \int_{\Omega}\rho(x)\nabla u(x)\cdot \nabla v(x)\,dx, \qquad (f,v) = \int_{\Omega}f(x) v(x)\,dx.
\]

We consider a family $\{ \Th \}_h$ of tessellations of $\Omega$ into a finite number of simple polygons $K$, and let $\Eh$ be the set of edges $e$ of $\Th$. For each tessellation $\Th$, we assume there exist constants $\gamma_0, \gamma_1, \alpha_0, \alpha_1 > 0$ such that: 
\begin{itemize}
\item each element $K \in \Th$ is star-shaped with respect to a ball of radius $\geq \gamma_0 h_K$, where $h_K$ is the diameter of $K$;
\item for each element $K\in\Th$ the distance between any two vertices of $K$ is $\geq \gamma_1 h_K$;
\item $\Th$ is quasi-uniform, that is, for any two elements $K$ and $K'$ in $\mathcal{T}_h$ we have $\alpha_0  \leq h_K / h_{K'}\leq \alpha_1 $.
\end{itemize}
%

For each polygon $K \in \Th$ we define a local finite element space $V^{K,k}$ as
\begin{equation*}
	V^{K,k} = \{ v \in H^1(K): v|_{\partial K} \in C^0(\partial K),\ v|_{e} \in \mathbb{P}_k(e)\ \forall e \in \EK,\ \Delta v \in \mathbb{P}_{k-2}(K) \},
\end{equation*}
with $\mathbb{P}_{-1}=\{0\}$. Then, the global virtual element space $V_h$ is defined as
\begin{align*}
V_h &= \{ v \in V: w|_{K} \in V^{K,k}\ \forall K \in \Th \}.
\end{align*}
We will consider  the following degrees of freedom, uniquely identifying a function $v_h \in V_h$:
\begin{itemize}
\item the values of $v_h$ at the vertices of the tessellation;
\item for each edge $e$, the values of $v_h$ at the $k-1$ internal points of the $k+1$ points Gauss-Lobatto quadrature rule on $e$;
\item for each element $K$, the moments up to order $k-2$ of $v_h$ in $K$. 
\end{itemize}
Due to the definition of the discrete space $V_h$, the bilinear form $a$ in equation~\eqref{eq:variational} is not directly computable on discrete functions in terms of the degrees of freedom. The VEM stems from replacing $a$ with a suitable approximate bilinear form $a_h$. Thus, the discrete form of problem~\ref{eq:poisson} reads as follows: find $u_h\in V_h$ such that
\begin{equation}\label{eq:poisson-h}
a_h(u_h,v_h) = f_h(v_h)\qquad\forall\,v_h\in V_h.
\end{equation}
Further details on how the bilinear form $a_h$, as well as the study of the convergence, stability and robustness properties of the method can be found in~\cite{basicVEM,beirao_var_coef,beirao_hp}.
For further details 
on the implementation of the method we refer to~\cite{beirao2014}.

\section{The FETI-DP Domain Decomposition Method for the VEM}\label{sec:feti-dp}

Since the degrees of freedom corresponding to the edges of the polygons in the tessellation are nodal values, the FETI-DP method is defined as in the finite element case. More precisely let  $\Omega$  be split as $\Omega = \cup_\ell \oml$, with $\oml = \cup_{K \in \Thk} K$, where $\Thk$ are disjoint subsets of $\Th$.  In view of the quasi uniformity assumptions on the tessellation $\Th$ and assuming that also the  decomposition into subdomains is quasi uniform, we can introduce global mesh size parameters $H$ and $h$ such that for all $\ell$ and for all $K$
we have 
$h_K \simeq h$ and $\text{diam}(\Omega^\ell) \simeq H$.
We let $\Gamma = \cup_\ell \partial\Omega^\ell \setminus \partial\Omega$ denote the skeleton (or interface) of the decomposition.  

Let $\widetilde V_h \supset V_h$ denote the space obtained by dropping the continuity constraint at all nodes interior to the macro edges of the decomposition (which we will call dual nodes), while retaining continuity at cross points. Problem~\ref{eq:poisson-h} is equivalent to finding $\tilde u_h\in\widetilde V_h$ satisfying
\begin{equation}\label{eq:min}
	\left\{
	\begin{array}{l}
		J(\tilde u_h) = \min_{\tilde v_h\in\widetilde V_h} J(\tilde v_h), \text{ with } J(\tilde v_h) = \dfrac{1}{2}a_h(\tilde v_h, \tilde v_h) - \int_\Omega f \tilde v_h,\\
		\text{such that $\tilde u_h$ is continuous across the interface.}
	\end{array}\right.
\end{equation}
For each $\tilde v_h\in\widetilde V_h$, we denote by $\vec{\tilde v}\in\mathbb R^M$ the corresponding vector coefficient, where $M$ is the dimension of $\widetilde V_h$. Let $B$ be a matrix whose entries assume value in the set $\{-1,0,1\}$. The continuity constraints across the interface can then be expressed as $B\vec{\tilde{u}} = \vec{0}$. By introducing a set of Lagrange multipliers $\boldsymbol\lambda \in \text{range}(B)$ to enforce the continuity constraints, we obtain a finite dimensional saddle point formulation of~\eqref{eq:min}
\[
\widetilde A\vec{\tilde{u}} + B^T\boldsymbol\lambda = \vec{\tilde f},\qquad B\vec{\tilde u} = \vec{0},
\]
where $\widetilde A$ and $\vec{\tilde f}$ are the finite dimensional representation of $a_h(\tilde v_h, \tilde w_h)$ and $\int_\Omega f\tilde v_h$, $\forall\,\tilde v_h, \tilde w_h\in\widetilde V_h$, respectively. Since $\widetilde A$ is symmetric and positive definite~\cite{toselli2005}, we can eliminate $\vec{\tilde u}$, and obtain a linear system for the Lagrange multiplier. This linear system is solved with a conjugate gradient method with a preconditioner that takes the same form as in the finite element case (\cite{toselli2005}). In \cite{Bertoluzza2017} the authors proved that, as in the finite element case, the condition number of the preconditioned matrix increases at most as $(1+\log(k^2H/h))^2$, under the assumption that the $\Omega^\ell$ are shape regular.

\section{Subdomain Partitioning by Conformal Meshing}\label{sec:conformal}

In general, subdomains $\Omega^\ell$ obtained as the union of polygonal elements of a tessellation are not shape regular, unless this is constructed in two stages: first, a decomposition into shape regular subdomains is defined; then, each subdomain is refined to obtain the final tessellation. An alternative is to define the subdomains independently of the tessellation and to modify the latter by ``cutting'' it with the edges of the subdomains. By construction, the subdomains will then be the union of elements of the new tessellation.

We provide here details of our meshing algorithm, which we implemented in C++ using CinoLib \cite{cinolib}. Given a tessellation $\mathcal{T}_h$ of $\Omega$ and a set of $L$ polygonal subdomains that jointly cover $\Omega$ without overlaps, we output a domain tessellation with matching interfaces along subdomain boundaries.

\begin{figure}
\centering
\includegraphics[width=\linewidth]{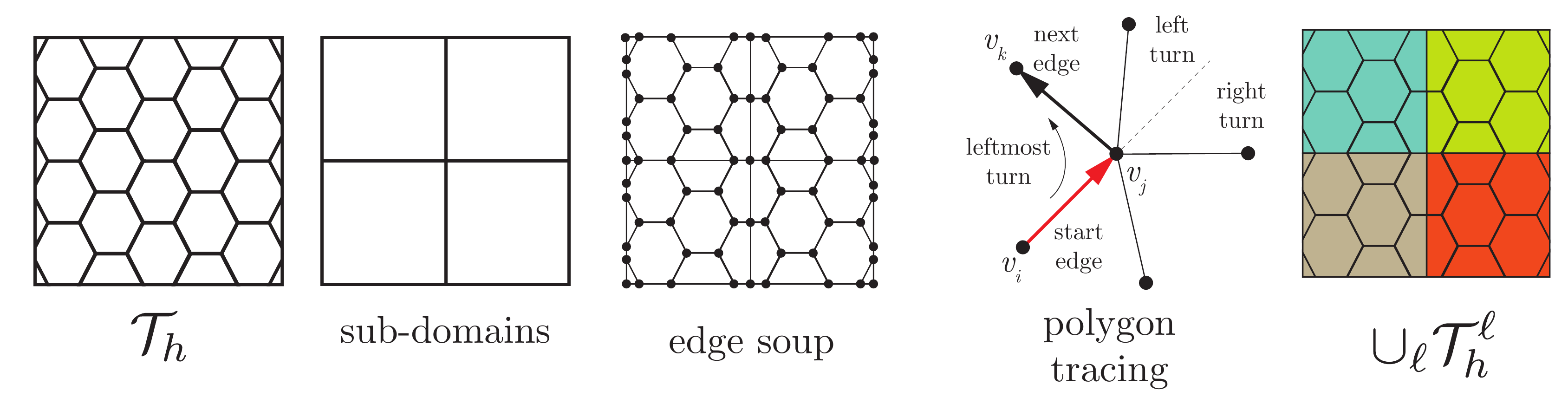}
\caption{Schematic representation of our meshing algorithm.}
\label{fig:meshing}
\end{figure}

In the general case the edges of $\mathcal{T}_h$ will intersect subdomain boundaries at many points. We start by splitting all the edges that are crossed by subdomain boundaries, thereby producing an edge soup. For any pair of edges in the soup, either the two edges are completely disjoint or they meet at a common endpoint (Fig.~\ref{fig:meshing}, middle). When computing intersections, exact arithmetics should be used to handle pathological cases due to roundoff errors introduced by finite machine arithmetics~\cite{repair}. 

In order to generate the output, we process the edge soup and produce a list of polygons, each one represented as an ordered list of vertices. To this end, for each oriented edge in the soup, we trace polygons by iteratively visiting the leftmost edge until hitting the starting edge at its opposite endpoint, thus closing the polygon (Fig.~\ref{fig:meshing}, middle right). If all the edges are visited twice (both orientations count) the result is a list of polygons, with vertices ordered counterclockwise. The only polygon with clockwise winding order will be the one covering the whole domain and containing all and only the boundary edges of $\mathcal{T}_h$, which is discarded during the post-processing phase.

The leftmost edge is found as follows: given an oriented edge $e$ (from node $v_i$ to node $v_j$), the set of edges incident to it at $v_j$ are first classified as being either on the left or on the right of $e$. Let $e'(v_j \rightarrow v_k)$ be one of these edges and $A$ the $2 \times 2$ matrix having $(v_j - v_i)$ and $(v_k - v_j)$ as rows. $e'$ is at the left of $e$ if $\det(A) > 0$, at the right of $e$ if $\det(A) < 0$ and collinear with $e$ if $\det(A) = 0$. Computing $\det(A)$ is known to be a critical issue in finite floating point arithmetics. For this reason, Shewchuck predicates \cite{shewchuk1997adaptive} are used to estimate its sign robustly. If there exist some edges at the left of $e$, then the polygon is locally convex and the leftmost turn is the one that passes through the edge $e'$ that minimizes  the dot product between $e' / \Vert e' \Vert $ and $e / \Vert e \Vert $. Conversely, if all the edges incident to $v_j$ are at the right of $e$, the polygon is locally concave and the leftmost turn is the one that passes through the edge $e''$ that maximizes the dot product between $e'' / \Vert e'' \Vert $ and $e / \Vert e \Vert $.

\section{Numerical Results}\label{sec:results}

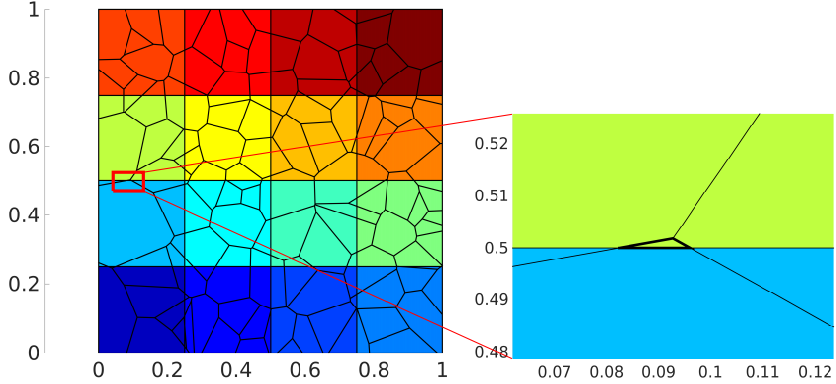
\begin{figure}
\centering
\begin{tikzpicture}
\node[immagine] at (0,0) {\includegraphics[width=0.5\textwidth]{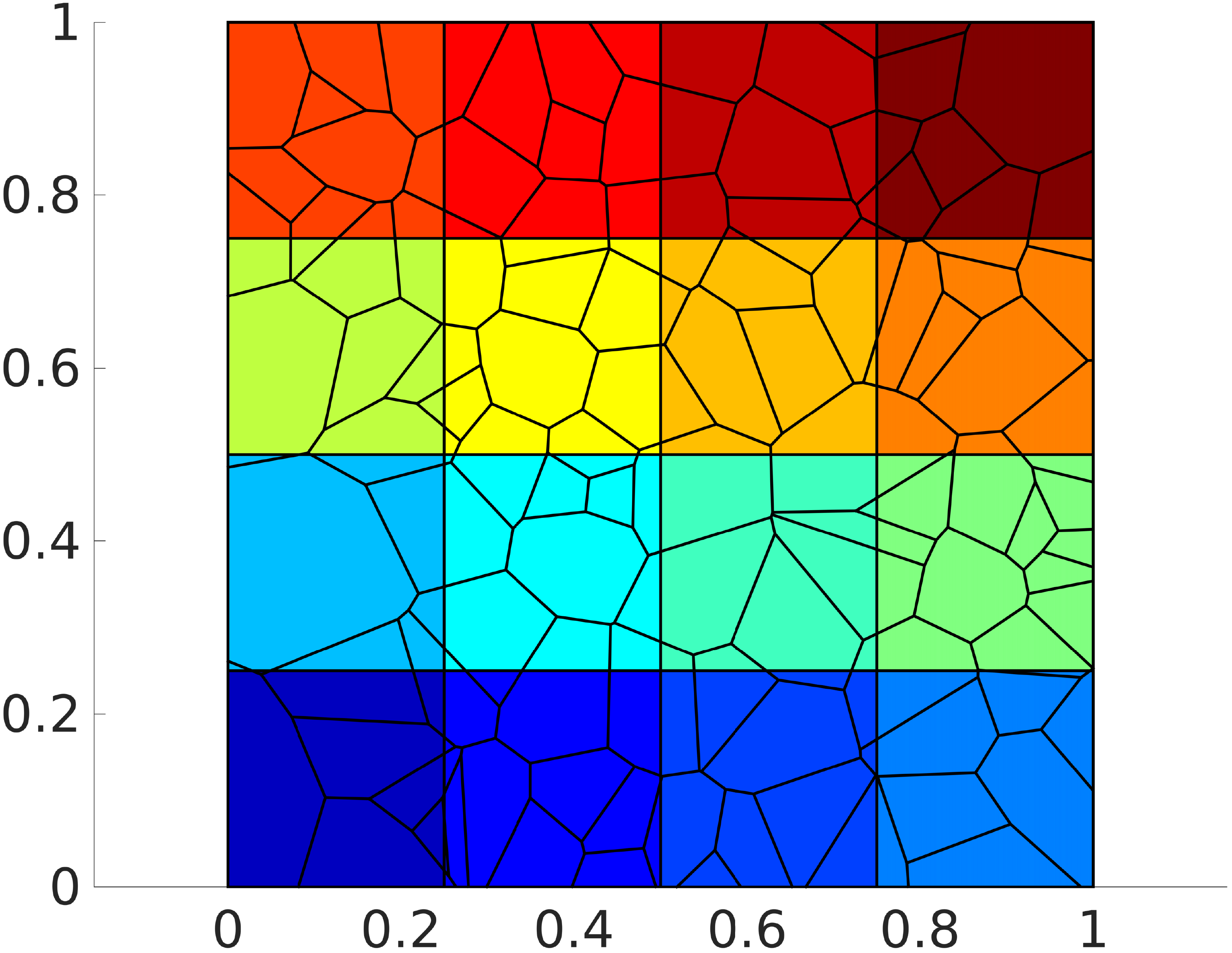}};
\node[immagine] at (6.2,0) {\includegraphics[width=4.8cm]{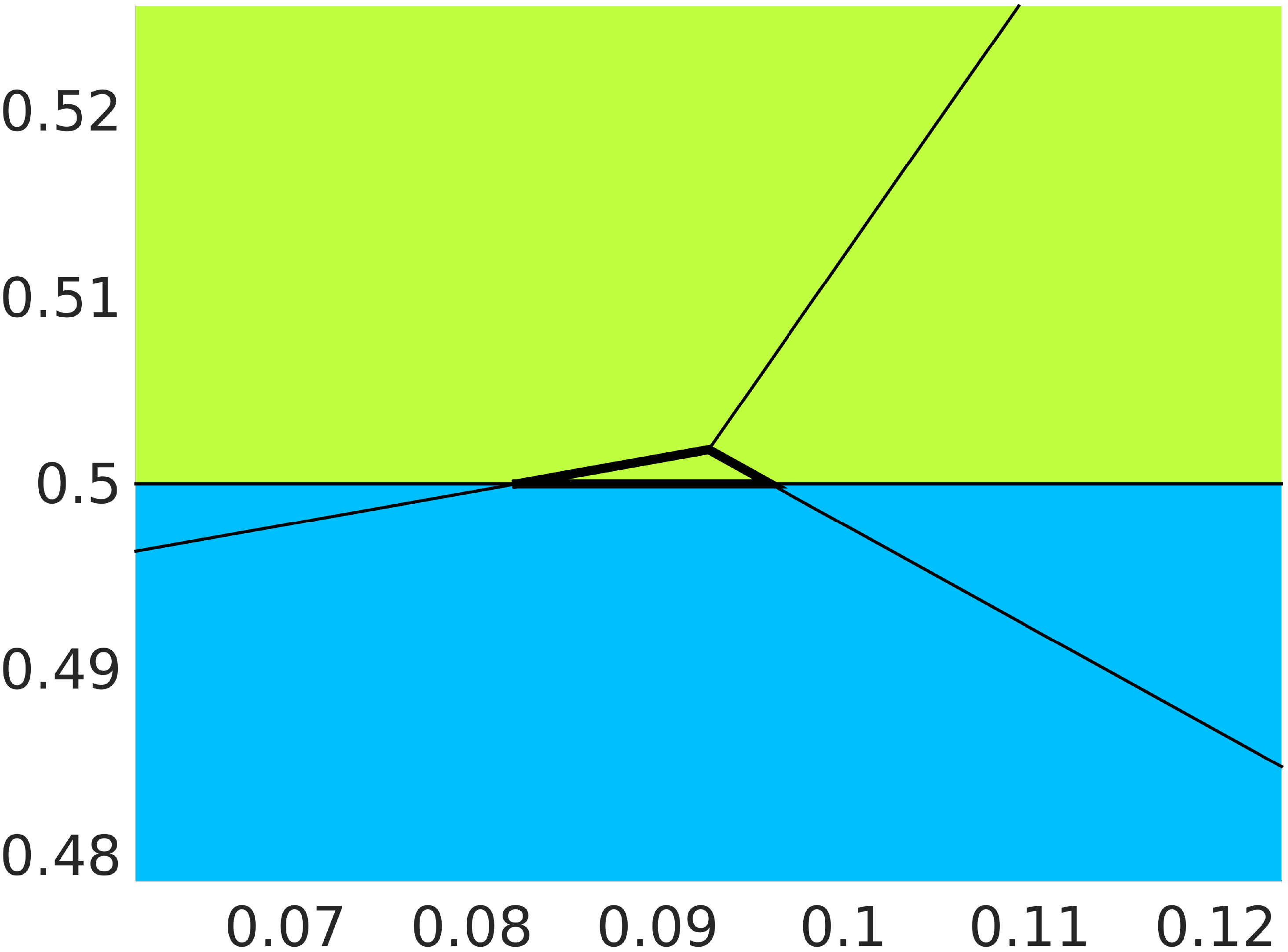}};
\draw[very thick, red] (1.4,2.5) -- (1.8,2.5) -- (1.8,2.75) -- (1.4,2.75) -- (1.4,2.5);
\draw[red] (1.8,2.5) -- (6.7,0.27);
\draw[red] (1.8,2.75) -- (6.7,3.52);
\end{tikzpicture}
\caption{(Left) General Voronoi mesh partitioned into $L = 16$ subdomains. (Right) Zoomed view of the boundary between two subdomains, where a tiny triangle is generated.}
\label{fig:tiny}
\end{figure}

Here we consider the performance of the preconditioned FETI-DP method. In each experiment, the domain $\Omega = [0,1]^2$ has been discretized using Voronoi cells of arbitrary shape. Every tessellation $\Th$ of $\Omega$ is partitioned into $L$ squared subdomains using the approach described in Section~\ref{sec:conformal}. This approach can introduce new polygons with very small edges, as shown in Fig.~\ref{fig:tiny}.
In order to test the robustness of FETI-DP, we consider two different types of data: i) $\rho = 1, f = \sin(2\pi x)\sin(2\pi y)$; ii) for each subdomain, $\rho = 10^\alpha, \alpha$ random integer in $[-5,5]$, $f$ uniform random in $[-1,1]$.

\begin{figure}
\centering
\begin{tabular}{c}


\begin{tikzpicture}
      \begin{semilogxaxis}[xlabel=Dof, legend pos=south east]
        \addplot [thick, color=red, mark=square, solid] coordinates {
        	(23257, 2.0958651642610966)
            (44548, 2.222975717036311)
            (86202, 2.346599760197499)
            (168223, 2.488401708028515)
            (330104, 2.63362866966488)
        };
        \addlegendentry{$L=64$}
        \addplot [thick, color=green, mark=*, solid] coordinates {
        	(25567, 2.062136974343745)
            (47801, 2.1162290113416593)
            (90561, 2.2619696820106445)
            (174136, 2.3878890042591747)
            (338451, 2.496599385412446)
        };
        \addlegendentry{$L=144$}
        \addplot [thick, color=blue, mark=x, solid] coordinates {
        	(27979, 1.924103291873768)
            (51000, 1.9825608538179087)
            (94954, 2.152855381631862)
            (180377, 2.3237514144083815)
            (346829, 2.4389422875159044)
        };
        \addlegendentry{$L=256$}
        \addplot [thick, color=red, mark=square, dashed] coordinates {
        	(23257, 1.9771123366426608)
            (44548, 2.2077301495585138)
            (86202, 2.237622072940415)
            (168223, 2.4229825220898267)
            (330104, 2.4616963632980107)
        };
        \addplot [thick, color=green, mark=*, dashed] coordinates {
        	(25567, 1.9953007629461585)
            (47801, 2.041478948384898)
            (90561, 2.2812433734178685)
            (174136, 2.3362653438062453)
            (338451, 2.4106170952273653)
        };
        \addplot [thick, color=blue, mark=x, dashed] coordinates {
        	(27979, 1.8633392242613998)
            (51000, 2.0208875308609664)
            (94954, 2.1371707785301717)
            (180377, 2.3231083145773064)
            (346829, 2.4151053886886356)
        };
      \end{semilogxaxis}
\end{tikzpicture}
\\
\end{tabular}
\caption{Plots of $\kappa^{1/2}$ as a function of the global degrees of freedom for fixed polynomial order $k = 1$ but increasing number of subdomains $L$. Solid and dashed lines correspond to the first and second types of data, respectively.}
\label{fig:kappa_rho1}
\end{figure}
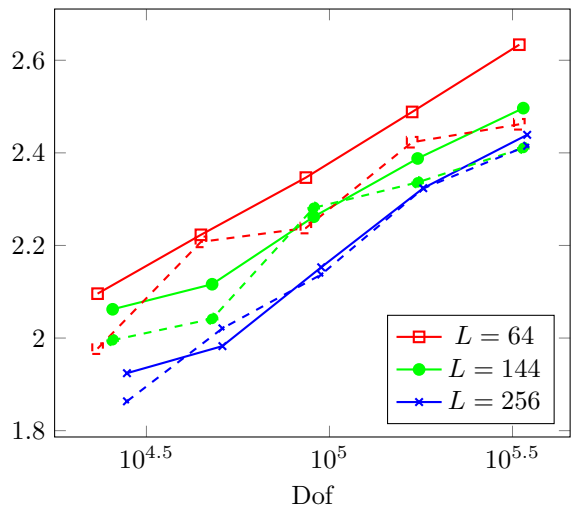

\begin{figure}
\centering
\begin{tabular}{rl}

\begin{tikzpicture}[trim axis left]
	\begin{semilogxaxis}[small, xlabel=$k$, ylabel=$\kappa^{1/2}$, legend pos=south east, xtick=data, xticklabels={$2$,$3$,$4$,$5$,$6$}]
    \addplot [thick, color=magenta, mark=*, solid] coordinates {
        (2, 3.1333590337331896e+00)
        (3, 3.5108207659958479e+00)
        (4, 3.7524602372942897e+00)
        (5, 3.9143480405395477e+00)
        (6, 4.0505956697079899e+00)
	};
    \addlegendentry{$L=64$}
    \addplot [thick, color=cyan, mark=square, solid] coordinates {
		(2, 2.9705441651111912e+00)
        (3, 3.3684034854866125e+00)
        (4, 3.5715477759193419e+00)
        (5, 3.7404732619026855e+00)
        (6, 3.8778374694296538e+00)
	};    
	\addlegendentry{$L=256$}
	\addplot [thick, color=magenta, mark=*, dashed] coordinates {
        (2, 2.9274604671036109e+00)
        (3, 3.2831101055858265e+00)
        (4, 3.5900078501389641e+00)
        (5, 3.7163646809944293e+00)
        (6, 3.8490366498705244e+00)
    };
	\addplot [thick, color=cyan, mark=square, dashed] coordinates {
    	(2, 3.0374600795372890e+00)
        (3, 3.4654489689831167e+00)
        (4, 3.7398310397969703e+00)
        (5, 3.8840354360818354e+00)
        (6, 4.0324714603164891e+00)
	};
	\end{semilogxaxis}
\end{tikzpicture}

&

\begin{tikzpicture}[trim axis right]
	\begin{semilogxaxis}[small, xlabel=Dof, legend pos=south east]
    	\addplot [thick, color=olive, mark=*, solid] coordinates {
        	(153660, 2.7325322891312487e+00)
            (280708, 2.8511293699631648e+00)
            (523710, 2.9926871202317380e+00)
            (996034, 3.1341323476392970e+00)
            (1916865, 3.3684034854866125e+00)
		};
        \addlegendentry{$k = 3$}
        \addplot [thick, color=teal, mark=square, solid] coordinates {
            (335433, 3.1708328392729412e+00)
            (612896, 3.2607463590910633e+00)
            (1143858, 3.3979356916824304e+00)
            (2175927, 3.5430565763474915e+00)
            (4188317, 3.7404732619026855e+00)
        };
        \addlegendentry{$k = 5$}
        \addplot [thick, color=olive, mark=*, dashed] coordinates {
            (153660, 3.0537005137087823e+00)
            (280708, 2.9024852539081523e+00)
            (523710, 3.2284941705884229e+00)
            (996034, 3.4229768190556182e+00)
            (1916865, 3.4654489689831167e+00)
        };
        \addplot [thick, color=teal, mark=square, dashed] coordinates {
            (335433, 3.5154207002174571e+00)
            (612896, 3.4744643061973881e+00)
            (1143858, 3.7282126848007522e+00)
            (2175927, 3.9081152195261026e+00)
            (4188317, 3.8840354360818354e+00)
        };        
	\end{semilogxaxis}
\end{tikzpicture}

\\
\end{tabular}
\caption{High order elements. Solid and dashed lines correspond to the first and second types of problem data, respectively. (Left) Plot of $\kappa^{1/2}$ as a function of the polynomial order $k$, initial mesh voro2 (see Table~\ref{tab:test1}). We have $H/h \approx 20.04$ for $L = 64$ and $H/h \approx 10.02$ for $L = 256$. (Right) Plot of $\kappa^{1/2}$ as a function of the global degrees of freedom, 
$H/h \approx 10.02, L = 256$.}
\label{fig:high-order}
\end{figure}
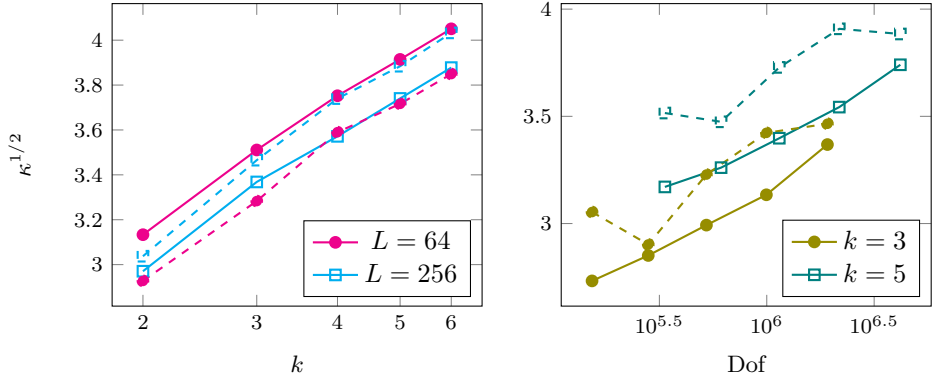

Table~\ref{tab:test1} shows that, with the first type of data, by fixing the initial mesh $\Th$ and increasing the number of subdomains, thereby fixing $h$ while decreasing $H$, the condition number $\kappa$ decreases as expected. The smallest edge in the mesh $h_\textup{min}$, as well as the parameters $\gamma_0$ and $\gamma_1$ introduced in Section~\ref{sec:vem} are also listed to stress that the presence of arbitrarily shaped Voronoi polygons does not hinder the convergence of the method. In Fig.~\ref{fig:kappa_rho1}, $\kappa^{1/2}$ is plotted as a function of the degrees of freedom of the whole problem, for varying number of subdomains $L$. The size of the problem is increased by taking five different meshes with $10000$, $20000$, $40000$, $80000$, and $160000$ polygons, respectively, as initial tessellations of $\Omega$. Solid and dashed lines correspond to the first and second types of data, respectively. With the first type, for fixed $L$, $\kappa^{1/2}$ clearly grows as $\log(\text{degrees of freedom})$, in agreement with the theoretical bound. This behavior does not seem to be affected by jumps in the diffusion coefficient, especially for high $L$ (see Fig.~\ref{fig:kappa_rho1}, right). Finally, we present some computations performed with high order elements (see Table~\ref{tab:high-order} and Fig.~\ref{fig:high-order}). In Fig.~\ref{fig:high-order} (Left), it is possible to verify the polylogarithmic dependence of $\kappa$ on the polynomial order $k$. Indeed, for $H/h$ fixed, the expected bound $(1+\log(k^2 H / h)^2 \sim (1 +log(k^2))^2 \sim(1 + 2\log(k))^2 $ is confirmed. In Fig.~\ref{fig:high-order} (Right), we keep the polynomial order fixed to $k = 3$ and $k = 5$ and increase the dimension of the problem by using the five meshes mentioned above. The expected bound $(1+log(k^2 H/h))^2 \sim(1+log( H/h))^2$ (for fixed $k$) is confirmed for both types of problem data.

All these experiments demonstrate that the performance of the preconditioned FETI-DP method is robust with respect to jumps in the diffusion coefficients and shape of the mesh elements.

\begin{table}
\caption{Results obtained with the first type of data and polynomial order $k = 1$ on two Voronoi meshes.}
\label{tab:test1}
  \centering
  \begin{tabular}{
      S[table-format=3.1]
      S[table-format=6.0]
      S[table-format=4.2]
      S[table-format=1.{\roundPrecision}e-1]
      S[table-format=1.{\roundPrecision}e-1]
      S[table-format=1.{\roundPrecision}e-1]
      S[table-format=1.1]
      S[table-format=1.2]
      S[table-format=2.0]
    }
    \toprule
        
        \multicolumn{9}{c}{{voro1, $40000$ initial polygons}}\\
        \midrule
            {$L$} & {D.o.f.} & {$1/h$} & {$h_\textup{min}$} & {$\gamma_0$} & {$\gamma_1$} & {$\lambda_\textup{min}$} & {$\lambda_\textup{max}$} & {It.}\\
            \midrule
            64 & 86202 & 55.93 & 1.04e-5 & 5.73e-3 & 1.32e-3 & 1.05 & 5.78 & 14\\
            144 & 90561 & 55.93 & 1.04e-5 & 3.97e-3 & 1.32e-3 & 1.05 & 5.37 & 14\\
            256 & 94954 & 55.93 & 1.04e-5 & 2.94e-3 & 1.13e-3 & 1.06 & 4.90 & 13\\
			\midrule
            
            \multicolumn{9}{c}{{voro2, $160000$ initial polygons}}\\
            \midrule
                {$L$} & {D.o.f.} & {$1/h$} & {$h_\textup{min}$} & {$\gamma_0$} & {$\gamma_1$} & {$\lambda_\textup{min}$} & {$\lambda_\textup{max}$} & {It.}\\
                \midrule
                64 & 330104 & 113.38 & 1.01e-5 & 4.11e-4 & 2.47e-3 & 1.06 & 7.32 & 16\\
                144 & 338451 & 113.38 & 1.04e-5 & 4.11e-4 & 2.10e-3 & 1.05 & 6.56 & 15\\
                256 & 346829 & 113.38 & 1.00e-5 & 4.11e-4 & 1.87e-3 & 1.05 & 6.25 & 15\\
                \bottomrule
  \end{tabular}
\end{table}

\begin{table}
\caption{$\kappa^{1/2}$ and number of iterations of the preconditioned FETI-DP system with the first type of data, for a fixed starting mesh (voro2, see Table~\ref{tab:test1}) but increasing number of subdomains $L$ and polynomial order $k$.}
\label{tab:high-order}
\centering
\begin{tabular}{
	S[table-format=3.0]
    S[table-format=1.2]
    S[table-format=2.0]
    S[table-format=1.2]
    S[table-format=2.0]
    S[table-format=1.2]
    S[table-format=2.0]
    S[table-format=1.2]
    S[table-format=2.0]
    S[table-format=1.2]
    S[table-format=2.0]
    }
\toprule
{$L \backslash k$} & \multicolumn{2}{c}{$2$} & \multicolumn{2}{c}{$3$} & \multicolumn{2}{c}{$4$} & \multicolumn{2}{c}{$5$} & \multicolumn{2}{c}{$6$}\\
\midrule
64 & 3.13 & 19 & 3.51 & 21 & 3.75 & 23 & 3.91 & 23 & 4.05 & 24\\
144 & 2.98 & 18 & 3.37 & 21 & 3.60 & 22 & 3.78 & 23 & 3.93 & 24\\
256 & 2.97 & 18 & 3.37 & 21 & 3.57 & 22 & 3.74 & 23 & 3.88 & 24\\
\bottomrule
\end{tabular}
\end{table}

\

\

\section*{Acknowledgments}
This paper has been realized in the framework of the ERC Project CHANGE, which has received funding from the European Research Council (ERC) under the European Union’s Horizon 2020 research and innovation programme (grant agreement No 694515). The authors would also like to thank the members of the Shapes and Semantics Modeling Group at IMATI-CNR for fruitful discussions on conformal meshing.

\bibliographystyle{siamplain}
\bibliography{references}

\end{document}